\documentclass[a4paper,12pt]{article}
\usepackage{longtable}
\title{Toric Fano varieties with divisorial contractions to curves}

\date{}

\renewcommand{\thefootnote}{\fnsymbol{footnote}}

\author{{\sc Hiroshi Sato}}

\newtheorem{Thm}{Theorem}[section]
\newtheorem{Prop}[Thm]{Proposition}

\newtheorem{Lem}[Thm]{Lemma}

\newtheorem{Def}[Thm]{Definition}
\newtheorem{Rem}[Thm]{Remark}

\newcommand{\proof}{Proof. \quad}
\newcommand{\qed}{\hfill q.e.d.}

\newcommand{\G}{\mathop{\rm G}\nolimits}

\newcommand{\PC}{\mathop{\rm PC}\nolimits}

\begin{document}
\maketitle

\renewcommand{\thefootnote}{}
\footnote{\hspace*{1.5em} {\em $2000$ Mathematics Subject Classification\/}.
Primary 14M25; 
Secondary 14E30, 14J45.}

%%%%%%%%%%%%%%%%%%%%%%%%%%%%%%%%%%%%
\begin{abstract}
In this paper, we obtain a complete classification of smooth toric Fano varieties equipped with extremal contractions which contract divisors to curves for any dimension. As an application, we obtain a complete classification of smooth projective toric varieties which can be equivariantly blown-up to Fano along curves.
\end{abstract}
%%%%%%%%%%%%%%%%%%%%%%%%%%%%%%%%%%%%

%%%%%%%%%%%%%%%%%%%%%%%%%%%%%%%%%%%%%%%%%%%%%
\section{Introduction}\label{intro}
%%%%%%%%%%%%%%%%%%%%%%%%%%%%1

\thispagestyle{empty}

\hspace{.5cm} Toric {\em Fano} $d$-fold $X$ is a smooth projective toric $d$-fold whose anti-canonical divisor $-K_X$ is ample. Toric Fano $d$-folds are classified for $d\leq 4$ (see Batyrev \cite{batyrev4}, Oda \cite{oda2}, Sato \cite{sato1} and Watanabe-Watanabe \cite{watanabe1}).

On the other hand, Bonavero \cite{bonavero1} classified toric Fano $d$-folds equipped with extremal contractions which contract divisors to {\em points} for any $d$. As a next step for this result, in this paper, we obtain a complete classification of toric Fano $d$-folds equipped with extremal contractions which contract divisors to {\em curves} for any $d$ (see Section \ref{table1}). Moreover, similarly as in Bonavero \cite{bonavero1}, we can classify smooth projective toric $d$-folds which can be equivariantly blown-up to Fano along curves for any $d$ (see Section \ref{table2}).

The author wishes to thank Professors Shihoko Ishii and Tatsuhiro Minagawa for advice and encouragement.

%%%%%%%%%%%%%%%%%%%%%%%%%%%%%%%%%%%%%%%%%%%%%%%
\section{Primitive collections and primitive relations}
%%%%%%%%%%%%%%%%%%%%%%%%%%%%%%%%%%%%%%%%%%%%%%%

\hspace{5mm} In this section, we review the concepts of primitive collections and primitive relations. They are very useful. See Batyrev \cite{batyrev3}, \cite{batyrev4}, Casagrande \cite{casagrande1}, \cite{casagrande2} and Sato \cite{sato1} more precisely. For fundamental properties of the toric geometry, see Fulton \cite{fulton1} and Oda \cite{oda2}.

\begin{Def}

{\rm
Let $X$ be a smooth complete toric $d$-fold, $\Sigma$ the corresponding fan in $N:={\bf Z}^d$ and $\G(\Sigma)\subset N$ the set of primitive generators of $1$-dimensional cones in $\Sigma$. A subset $P\subset\G(\Sigma)$ is called a {\em primitive collection} of $\Sigma$ if $P$ does not generate a cone in $\Sigma$, while any proper subset of $P$ generates a cone in $\Sigma$. We denote by $\PC(\Sigma)$ the set of primitive collections of $\Sigma$.
}

\end{Def}

Let $P=\{x_1,\ldots,x_m\}$ be a primitive collection of $\Sigma$. Then, there exists a unique cone $\sigma(P)$ in $\Sigma$ such that $x_1+\cdots+x_m$ is contained in the relative interior of $\sigma(P)$, because $X$ is complete. So, we get an equality
$$x_1+\cdots+x_m=a_1y_1+\cdots+a_ny_n,$$
where $y_1,\ldots,y_n$ are the generators of $\sigma(P)$, that is, $\sigma(P)\cap\G(\Sigma)=\{y_1,\ldots,y_n\}$, and $a_1,\ldots,a_n$ are positive integers. We call this equality the {\em primitive relation} of $P$. Thus, we obtain an element $r(P)$ in $A_{1}(X)$ for any primitive collection $P\in\PC(\Sigma)$, where $A_{1}(X)$ is the group of $1$-cycles on $X$ modulo rational equivalences. We define the {\em degree} of $P$ as $\deg P:=\left( -K_X\cdot r(P)\right)=m-(a_1+\cdots+a_n)$. The following is important.

\begin{Prop}[Batyrev \cite{batyrev3}, Reid \cite{reid1}]\label{toriccone}

Let $X$ be a smooth projective toric variety and $\Sigma$ the corresponding fan. Then
$${\bf NE}(X)=\sum_{P\in\PC(\Sigma)}{\bf R}_{\geq 0}r(P),$$
where ${\bf NE}(X)$ is the Mori cone of $X$.

\end{Prop}

A primitive collection $P$ is called an {\em extremal} primitive collection when $r(P)$ is contained in an extremal ray of ${\bf NE}(X)$. In particular, if $\deg P=1$ for a primitive collection $P\in\PC(\Sigma)$, then $P$ is an extremal primitive collection.

\begin{Def}

{\rm
Let $X$ be a smooth projective algebraic $d$-fold. Then, $X$ is called a {\em Fano} $d$-fold, if its anti-canonical divisor $-K_X$ is ample.
}

\end{Def}

By using the notion of primitive collections and primitive relations, toric Fano $d$-folds are characterized as follows.

\begin{Prop}[Batyrev \cite{batyrev4}, Sato \cite{sato1}]

Let $X$ be a smooth projective toric $d$-fold and $\Sigma$ the corresponding fan. Then, $X$ is Fano if and only if $\deg P>0$ for any primitive collection $P\in\PC(\Sigma)$.

\end{Prop}

We review some important results about toric Fano $d$-folds, primitive collections and primitive relations in Casagrande \cite{casagrande1}, \cite{casagrande2} and Sato \cite{sato1}. They are necessary for the classification.

\begin{Prop}[Casagrande \cite{casagrande1}, Sato \cite{sato1}]\label{contract}

Let $X$ be a smooth projective toric $d$-fold, $\Sigma$ the corresponding fan and $P$ a primitive collection of $\Sigma$ with primitive relation $x_1+\cdots+x_m=a_1y_1+\cdots+a_ny_n$, where $\{ x_1,\ldots,x_m,y_1,\ldots,y_n\}\subset\G(\Sigma)$ and $a_1,\ldots,a_n\in{\bf Z}_{>0}$. If $P$ is extremal, then, for any $P'\in\PC(\Sigma)$ such that $P\cap P'\neq\emptyset$ and $P\neq P'$, the set $(P'\setminus P)\cup\{y_1,\ldots,y_n\}$ contains a primitive collection.

\end{Prop}

\begin{Prop}[Casagrande \cite{casagrande2}]\label{redpri}

Let $X$ be a toric Fano $d$-fold and $\Sigma$ the corresponding fan. Suppose that $\{x, (-x)\}\in\PC(\Sigma)$. Then, for a primitive collection $P=\{x,y_1,\ldots,y_m\}\in\PC(\Sigma)$ containing $x$ such that $P\neq\{x,(-x)\}$, its primitive relation is
$$x+y_1+\cdots+y_m=z_1+\cdots+z_m,$$
where $\{z_1,\ldots,z_m\}\subset\G(\Sigma)$. Moreover, $\{(-x),z_1,\ldots,z_m\}$ is also a primitive collection of $\Sigma$ and its primitive relation is
$$(-x)+z_1+\cdots+z_m=y_1+\cdots+y_m.$$
These primitive relations are extremal.

\end{Prop}

\begin{Prop}[Casagrande \cite{casagrande2}]\label{pri2}

Let $X$ be a toric Fano $d$-fold and $\Sigma$ the corresponding fan. If there exist two distinct primitive relations $x+y=z$ and $x+w=v$, where $\{x,y,z,v,w\}\subset\G(\Sigma)$, then $w=(-z)$ and $v=(-y)$. In particular, we have the primitive relations
$$x+y=z,\ x+(-z)=(-y),\ y+(-y)=0,\ z+(-z)=0\mbox{ and }z+(-y)=x.$$

\end{Prop}

\begin{Thm}[Casagrande \cite{casagrande2}]\label{casfano}

Let $X$ be a toric Fano $d$-fold and $D$ a toric prime divisor on $X$. Then, we have $0\leq\rho(X)-\rho(D)\leq 3$, where $\rho(X)$ $($resp. $\rho(D))$ is the Picard number of $X$ $($resp. $D)$. Moreover, if $\rho(X)-\rho(D)=3$, then $X$ is an $S_6$-bundle over a toric Fano $(d-2)$-fold, where $S_6$ is the del Pezzo surface of degree $6$.

\end{Thm}

We close this section by giving the following fundamental definition and proposition.

\begin{Def}

{\rm
Let $X$ be a smooth complete toric $d$-fold and $\Sigma$ the corresponding fan. $\Sigma$ is called a {\em splitting fan} if $P\cap P'=\emptyset$ for any distinct primitive collections $P,P'\in\PC(\Sigma)$
}

\end{Def}

\begin{Prop}[Batyrev \cite{batyrev3}]

Let $X$ be a smooth complete toric $d$-fold and $\Sigma$ the corresponding fan. Then, $\Sigma$ is a splitting fan if and only if there exists a sequence of smooth complete toric varieties
$$X=X_1\stackrel{\psi_1}{\rightarrow}X_2\stackrel{\psi_2}{\rightarrow}\cdots\stackrel{\psi_{r-1}}{\rightarrow}X_{r}\stackrel{\psi_{r}}{\rightarrow}X_{r+1}$$
such that $\psi_{i}$ is a toric projective space bundle structure for $1\leq i\leq r$, while $X_{r+1}$ is a projective space.

\end{Prop}

%%%%%%%%%%%%%%%%%%%%%%%%%%%%%%%%%%%%%%%%%%%%%%%
\section{Classification}
%%%%%%%%%%%%%%%%%%%%%%%%%%%%%%%%%%%%%%%%%%%%%%%

\hspace{.5cm} In this section, we give the classification of toric Fano $d$-folds equipped with extremal contractions which contract divisors to curves. Since toric Fano $d$-folds are classified for $d\leq 4$, we assume $d\geq 5$ throughout this paper.

Let $X$ be a toric Fano $d$-fold, $\Sigma$ the corresponding fan and $\varphi:X\rightarrow Y$ an extremal contraction which contracts a divisor $E$ to a curve. We use this notation throughout this section. Since $E$ is a toric prime divisor on $X$, we have
$$0\leq\rho(X)-\rho(E)\leq 3$$
by Theorem \ref{casfano}. On the other hand, since $\varphi(E)$ is isomorphic to ${\bf P}^1$, we have $\rho(E)=2$. Therefore, the above inequalities are
$$2\leq\rho(X)\leq 5.$$
So, we consider the classification for these four cases separately. Let $x_1+\cdots+x_{d-1}=\alpha x$ be the extremal primitive relation corresponding to $\varphi$, where $\{ x_1,\cdots,x_{d-1},x\}\subset\G(\Sigma)$ and $1\leq\alpha\leq d-2$, $\{y_1,y_2\}\subset\G(\Sigma)$ the distinct elements such that $\{x_1,\ldots,x_{d-2},x,y_1\}$ and $\{x_1,\ldots,x_{d-2},x,y_2\}$ generate maximal cones of $\Sigma$, and $\{z_1,\ldots,z_{\rho(X)-2}\}\subset\G(\Sigma)$ the other elements.

\bigskip

(I) \underline{$\rho(X)=2$.}

\bigskip

Kleinschmidt \cite{klein1} showed that the fan of a smooth complete toric $d$-fold of Picard number $2$ is a splitting fan. Therefore, in this case, $X$ is a ${\bf P}^{2}$-bundle over a ${\bf P}^{d-2}$. So, we can easily determine the corresponding fan (see Section \ref{table1}).

\bigskip

(II) \underline{$\rho(X)=3$.}

\bigskip

Smooth projective toric $d$-folds of Picard number $3$ are classified in Batyrev \cite{batyrev3}. There are explicit descriptions of fans as follows.

\begin{Thm}[Batyrev \cite{batyrev3}]\label{baty3}

Let $X$ be a smooth projective toric $d$-fold of Picard number $3$ and $\Sigma$ the corresponding fan. Then, one of the following holds.
\begin{enumerate}
\item $\Sigma$ is a splitting fan.
\item $\#\PC(\Sigma)=5$.
\end{enumerate}
Moreover, in the case of $(2)$, there exists $(p_0,p_1,p_2,p_3,p_4)\in({\bf Z}_{>0})^{5}$ such that the primitive relations of $\Sigma$ are
$$v_1+\cdots+v_{p_0}+s_{1}+\cdots+s_{p_1}=c_2w_2+\cdots+c_{p_2}w_{p_2}+(b_1+1)t_1+\cdots+(b_{p_3}+1)t_{p_3},$$
$$s_1+\cdots+s_{p_1}+w_1+\cdots+w_{p_2}=u_1+\cdots+u_{p_4},\ w_1+\cdots+w_{p_2}+t_1+\cdots+t_{p_3}=0,$$
$$t_1+\cdots+t_{p_3}+u_1+\cdots+u_{p_4}=s_{1}+\cdots+s_{p_1}\mbox{ and}$$
$$u_1+\cdots+u_{p_4}+v_1+\cdots+v_{p_0}=c_2w_2+\cdots+c_{p_2}w_{p_2}+b_1t_1+\cdots+b_{p_3}t_{p_3},$$
where $\G(\Sigma)=\{v_1,\ldots,v_{p_0},s_{1},\ldots,s_{p_1},w_1,\ldots,w_{p_2},t_1,\ldots,t_{p_3},u_1,\cdots,u_{p_4}\}$ and $c_2,\ldots,c_{p_2},$ $b_1,\ldots,b_{p_3}\in{\bf Z}_{\geq 0}.$

\end{Thm}

Using this theorem, we can determine fans completely (see Section \ref{table1}).

\bigskip

%%%%%%%%%%%%%%%%%%%%%%%%%%%%%%%%%%%%%%%%%%%%%%%

(III) \underline{$\rho(X)=4$.}

\bigskip

This case is the main part of the classification. First, we need the following proposition.

\begin{Prop}\label{dp7}

Let $X$ be a toric Fano $d$-fold and $\Sigma$ the corresponding fan. Suppose that $\G(\Sigma)=\{x_1,\ldots,x_{d-1},u_1,u_2,u_3,v_1,v_2\}$, and that we have the primitive relations
$$u_1+u_3=u_2,\ u_1+v_1=0,\ u_2+v_1=u_3,\ u_2+v_2=0\mbox{ and }u_3+v_2=v_1,$$
and the extremal primitive relation
$$x_1+\cdots+x_{d-1}=\alpha x,$$
where $x\in\{u_1,u_2,u_3\}$ and $1\leq\alpha\leq d-2$. Then, $X$ is an $S_7$-bundle over ${\bf P}^{d-2}$, where $S_7$ is the del Pezzo surface of degree $7$.

\end{Prop}

To prove Proposition \ref{dp7}, we need the following lemmas.

\begin{Lem}\label{case4}

Let $P\in\PC(\Sigma)$ be a primitive collection. Then, the following hold.
\begin{enumerate}
\item If $P\cap\{x_1,\ldots,x_{d-1}\}\neq\emptyset$, then $x\not\in P$.
\item If $P\cap\{u_1,u_3\}\neq\emptyset$, then $u_2\not\in P$.
\item If $P\cap\{u_2,v_1\}\neq\emptyset$, then $u_3\not\in P$.
\item If $P\cap\{u_3,v_2\}\neq\emptyset$, then $v_1\not\in P$.
\end{enumerate}

\end{Lem}

\proof
We prove the case of $(1)$. The other cases are similar.

Suppose that $x\in P$. Since $x_1+\cdots+x_{d-1}=\alpha x$ is an extremal primitive relation, $(P\setminus\{x_1,\ldots,x_{d-1}\})\cup\{x\}$ contains a primitive collection by Proposition \ref{contract}. However, this is impossible, because $(P\setminus\{x_1,\ldots,x_{d-1}\})\cup\{x\}$ is a proper subset of the primitive collection $P$. So, we have $x\not\in P$.\qed

\begin{Lem}\label{exw}

Let $P\in\PC(\Sigma)$. If $P\cap\{x_1,\ldots,x_{d-1}\}\neq\emptyset$ and $P\cap\{u_1,u_2,u_3,v_1,v_2\}\neq\emptyset$, then there exists $w\in P\setminus\{x_1,\ldots,x_{d-1}\}$ such that $\{x,w\}$ is a primitive collection.

\end{Lem}

\proof
By Proposition \ref{contract}, $(P\setminus\{x_1,\ldots,x_{d-1}\})\cup\{x\}$ contains a primitive collection. On the other hand, obviously, there does not exist a primitive collection contained in $\{u_1,u_2,u_3,v_1,v_2\}$ other than $\{u_1,u_3\},\ \{u_1,v_1\},\ \{u_2,v_1\},\ \{u_2,v_2\}\mbox{ and }\{u_3,v_2\}.$\qed

\bigskip

Proof of Proposition \ref{dp7}. \quad Suppose that there exists a primitive collection $P\in\PC(\Sigma)$ such that $P\cap\{x_1,\ldots,x_{d-1}\}\neq\emptyset$ and $P\cap\{u_1,u_2,u_3,v_1,v_2\}\neq\emptyset$. Put $S=P\setminus\{x_1,\ldots,x_{d-1}\}$. We may assume $P=S\cup\{x_1,\ldots,x_l\}$ for $1\leq l<d-1$. We have to consider the following three cases.

\bigskip

$\underline{x=u_1.}$ By Lemma \ref{exw}, $u_3\in P$ or $v_1\in P$. For the case $u_3\in P$, we have $u_2\not\in P$ and $v_1\not\in P$ by $(2)$ and $(4)$ in Lemma \ref{case4}. Moreover, we have $v_2\not\in P$, since $S$ does not contain any primitive collection. So, we have $S=\{u_3\}$. Similarly, we have $S=\{v_1\}$ for the case $v_1\in P$.

Suppose that $S=\{u_3\}$. Since $u_1+u_3=u_2$ is an extremal primitive relation, $\{u_2,x_1,\ldots,x_l\}$ contains a primitive collection by Proposition \ref{contract}. This contradicts Lemma \ref{exw}.

Suppose that $S=\{v_1\}$. Since $u_2+v_1=u_3$ is an extremal primitive relation, $\{u_3,x_1,\ldots,x_l\}$ contains a primitive collection by Proposition \ref{contract}. This is impossible as above.

\bigskip

$\underline{x=u_2.}$ Similarly as in the case $x=u_1$, there exist three possibilities $S=\{v_1\}$, $S=\{v_2\}$ and $S=\{u_1,v_2\}$.

Suppose that $S=\{v_1\}$ or $S=\{v_2\}$. These cases are impossible by the similar argument as in the case $x=u_1$.

Let $S=\{u_1,v_2\}$. Since $u_2+v_2=0$, there exist two primitive relations
$$v_2+u_1+x_1+\cdots+x_l=w_1+\cdots+w_{l+1}\mbox{ and }u_2+w_1+\cdots+w_{l+1}=u_1+x_1+\cdots+x_{l}$$
by Proposition \ref{redpri}, where $\{w_1,\ldots,w_{l+1}\}\subset\G(\Sigma)$. Since $\#\{u_2,w_1,\ldots,w_{l+1}\}\geq 3$, we have $\{w_1,\ldots,w_{l+1}\}\cap\{x_1,\ldots,x_{d-1}\}\neq\emptyset$. This contradicts $(1)$ in Lemma \ref{case4}.

\bigskip

$\underline{x=u_3.}$ We can prove that this case is also impossible by the similar argument as in the case $x=u_2$.

\bigskip

Therefore, there does not exist a primitive collection other than
$$\{u_1,u_3\},\ \{u_1,v_1\},\ \{u_2,v_1\},\ \{u_2,v_2\},\ \{u_3,v_2\}\mbox{ and }\{x_1,\ldots,x_{d-1}\}.$$
Thus, $X$ is an $S_7$-bundle over ${\bf P}^{d-2}$.\qed

\bigskip

Now, we return to the classification. We have $\G(\Sigma)=\{x_1,\ldots,x_{d-1},x,y_1,y_2,z_1,z_2\}$. $P_1=\{x,z_1\}$ and $P_2=\{x,z_2\}$ are primitive collections. It is sufficient to consider the classification for the following four cases.
\begin{enumerate}
\item $\deg P_1=\deg P_2=1$.
\item $x+z_1=0$ and $x+z_2=y_1$.
\item $x+z_1=0$ and $x+z_2=z_1$.
\item $x+z_1=0$ and $x+z_2=x_1$.
\end{enumerate}

\bigskip

(1) \underline{$\deg P_1=\deg P_2=1$.} By Proposition \ref{pri2}, $\{(-z_1),(-z_2)\}\subset\G(\Sigma)$ and the primitive relations corresponding to $P_1$ and $P_2$ are
$$x+z_1=(-z_2)\mbox{ and }x+z_2=(-z_1),$$
respectively. If $(-z_1)=x_i$ for some $1\leq i\leq d-1$, then we have the primitive relation $x_i+z_1=0$. However, since $x_1+\cdots+x_{d-1}=\alpha x$ is an extremal primitive relation, this contradicts Proposition \ref{redpri} as long as $d\geq 4$. Therefore, $(-z_1)\in\{y_1,y_2\}$. Similarly, $(-z_2)\in\{y_1,y_2\}$. Let $y_1=(-z_1)$ and $y_2=(-z_2)$. Then, we have the primitive relations
$$y_1+z_1=0,\ y_2+z_2=0,\ x+z_1=y_2,\ x+z_2=y_1\mbox{ and }y_1+y_2=x.$$
So, $X$ is an $S_7$-bundle over ${\bf P}^{d-2}$ by Proposition \ref{dp7}.

\bigskip

(2) \underline{$x+z_1=0$ and $x+z_2=y_1$.} Obviously, $\{x,y_1,y_2\}$ contains a primitive collection. Suppose that $\{x,y_1,y_2\}$ is a primitive collection. Since $x+z_2=y_1$ is an extremal primitive relation, $\{y_1,y_2\}$ is a primitive collection by Proposition \ref{contract}. This is a contradiction. So, $\{y_1,y_2\}$ is a primitive collection.

$y_1+z_1=z_2$ is also an extremal primitive relation. So, $\{y_2,z_2\}$ is a primitive collection by Proposition \ref{contract}. Suppose that $y_1+y_2\neq 0$ and $y_2+z_2\neq 0$. Then, $(-y_1),\ (-z_2)\in\G(\Sigma)$, and hence we have the primitive relations $y_1+y_2=(-z_2)$ and $y_2+z_2=(-y_1)$ by Proposition \ref{pri2}. Therefore, $\{y_1,y_2\}$, $\{y_1,z_1\}$ and $\{y_1,(-y_1)\}$ are distinct primitive collections. So, we have $\rho(X)-\rho(D_1)=3$, where $D_1$ is the toric prime divisor corresponding to $y_1$. By Theorem \ref{casfano}, $X$ has an $S_6$-bundle structure. This is impossible.

If $y_1+y_2=0$, then we have the primitive relations
$$x+z_1=0,\ x+z_2=y_1,\ y_1+y_2=0,\ y_1+z_1=z_2\mbox{ and }y_2+z_2=z_1.$$
By Proposition \ref{dp7}, $X$ is an $S_7$-bundle over ${\bf P}^{d-2}$.

If $y_2+z_2=0$, then we have the primitive relations
$$x+z_1=0,\ x+z_2=y_1,\ y_1+y_2=x,\ y_1+z_1=z_2\mbox{ and }y_2+z_2=0.$$
By Proposition \ref{dp7}, $X$ is an $S_7$-bundle over ${\bf P}^{d-2}$.

\bigskip

(3) \underline{$x+z_1=0$ and $x+z_2=z_1$.} Since $x+z_2=z_1$ is an extremal primitive relation, we have a contradiction by Proposition \ref{contract}. So, this case is impossible.

\bigskip

(4) \underline{$x+z_1=0$ and $x+z_2=x_1$.} $\{x,y_1,y_2\}$ contains a primitive collection. Suppose that $\{x,y_1,y_2\}$ is a primitive collection. By Proposition \ref{redpri}, we have two extremal primitive relations
$$x+y_1+y_2=w_1+w_2\mbox{ and }z_1+w_1+w_2=y_1+y_2,$$
where $\{w_1,w_2\}\subset\G(\Sigma)$ and $w_1\neq w_2$. By Proposition \ref{contract}, we have $z_2\not\in\{w_1,w_2\}$, because $\{x,z_2\}$ is a primitive collection. Therefore, $\{w_1,w_2\}\subset\{x_1,\ldots,x_{d-1}\}$. Let $w_1=x_i$ and $w_2=x_j$ $(1\leq i<j\leq d-1)$. Since $x+z_2=x_1$ is an extremal primitive collection, $\{x_1,y_1,y_2\}$ contains a primitive collection by Proposition \ref{contract}. So, $\{x_1,y_1,y_2\}$ is a primitive collection. Since $x+y_1+y_2=x_i+x_j$ is an extremal primitive collection, $\{x_1,x_i,x_j\}$ contains a primitive collection by Proposition \ref{contract}. However, this is impossible as long as $d\geq 5$. So, $\{y_1,y_2\}$ is a primitive collection.

$x_1+z_1=z_2$ is an extremal primitive relation. So, $\{x_2,\ldots,x_{d-1},z_2\}$ is also a primitive collection by Proposition \ref{contract}. Thus, we have the primitive relations
$$x_1+\cdots+x_{d-1}=\alpha x,\ x_2+\cdots+x_{d-1}+z_2=(\alpha-1)x,$$
$$x+z_1=0,\ x+z_2=x_1\mbox{ and }z_1+x_1=z_2.$$
Therefore, $\Sigma$ contains a subfan $\Sigma'$ such that $\G(\Sigma')=\{x_1,\ldots,x_{d-1},x,z_1,z_2\}$ and the corresponding toric $(d-1)$-fold $X'$ is a toric Fano $(d-1)$-fold equipped with an extremal contraction which contracts a divisor to a point. Since $\{y_1,y_2\}$ is a primitive collection, $X$ is an $X'$-bundle over ${\bf P}^1$.

\bigskip

Thus, we obtain the following theorem.

\begin{Thm}\label{r4}

Let $X$ be a toric Fano $d$-fold of Picard number $4$. If there exists an extremal contraction from $X$ which contracts a divisor to a curve, then $X$ is either an $S_7$-bundle over ${\bf P}^{d-2}$ or an $X'$-bundle over ${\bf P}^1$, where $X'$ is a toric Fano $(d-1)$-fold equipped with an extremal contraction which contracts a divisor to a point.

\end{Thm}

By Theorem \ref{r4}, we can describe $\Sigma$ explicitly (see Section \ref{table1}).

%%%%%%%%%%%%%%%%%%%%%%%%%%%%%%%%%%%%%%%%%%%%%%%%

\bigskip

(IV) \underline{$\rho(X)=5$.}

\bigskip

By Theorem \ref{casfano}, $X$ is an $S_6$-bundle over a ${\bf P}^{d-2}$. So, we can easily determine the corresponding fan (see Section \ref{table1}).

%%%%%%%%%%%%%%%%%%%%%%%%%%%%%%%%%%%%%%%%%%%%%%%
\section{The classified list}\label{table1}
%%%%%%%%%%%%%%%%%%%%%%%%%%%%%%%%%%%%%%%%%%%%%%%

\hspace{.5cm} In this section, we give the complete list of toric Fano $d$-folds equipped with extremal contractions which contract divisors to curves. We assume $d\geq 5$. Let $1\leq \alpha\leq d-2$.

\bigskip

(I) \underline{$\rho(X)=2$.} $X$ is a ${\bf P}^2$-bundle ${\bf P}_{{\bf P}^{d-2}}\left(\mathcal{O}\oplus\mathcal{O}\oplus\mathcal{O}(\alpha)\right)$ over ${\bf P}^{d-2}$. The primitive relations are
$$x_1+\cdots+x_{d-1}={\alpha}x_d\mbox{ and }x_d+x_{d+1}+x_{d+2}=0,$$
where $\G(\Sigma)=\{x_1,\ldots,x_{d+2}\}$.

\bigskip

(IIa) \underline{$\rho(X)=3$ and $\Sigma$ is a splitting fan.} Let $\G(\Sigma)=\{x_1,\ldots,x_{d+3}\}$. The primitive relations are $x_{d+2}+x_{d+3}=0$ and
\begin{center}
\begin{tabular}{|c||c|c|c|c|c|} \hline
Case & 1 & 2 & 3 & 4 & 5 \\ \hline\hline
$x_1+\cdots+x_{d-1}=$ & ${\alpha}x_d$ & ${\alpha}x_d$ & ${\alpha}x_{d+2}$ & ${\alpha}x_{d+3}$ & $\alpha x_{d+2}$ \\ \hline
$x_d+x_{d+1}=$ & $0$ & $x_{d+2}$ & $x_{d+2}$ & $x_{d+2}$ & $x_1$ \\ \hline
\end{tabular}
\end{center}

\bigskip

(IIb) \underline{$\rho(X)=3$ and $\Sigma$ is not a splitting fan.} Put $\G(\Sigma)=\{x_1,\ldots,x_{d+3}\}$. There exist the following four cases.

\bigskip

(1) The primitive relations are $x_1+\cdots+x_{d-1}=\alpha x_d,\ x_1+\cdots+x_{d-3}+x_{d+1}+x_{d+2}=(\alpha-1)x_d,\ x_{d-2}+x_{d-1}+x_{d+3}=x_{d+1}+x_{d+2},\ x_d+x_{d+1}+x_{d+2}=x_{d-2}+x_{d-1}\mbox{ and }x_d+x_{d+3}=0$.

\bigskip

(2) The primitive relations are $x_1+\cdots+x_{d-1}=\alpha x_d,\ x_1+\cdots+x_{d-2}+x_{d+3}=(\alpha-1)x_d,\ x_{d-1}+x_{d+1}+x_{d+2}=x_{d+3},\ x_d+x_{d+1}+x_{d+2}=0\mbox{ and }x_d+x_{d+3}=x_{d-1}$.

\bigskip

(3) $\alpha=1$. The primitive relations are $x_1+\cdots+x_{d-1}=x_d,\ x_3+\cdots+x_{d-1}+x_{d+1}+x_{d+2}=0,\ x_1+x_2+x_{d+3}=x_{d+1}+x_{d+2},\ x_d+x_{d+1}+x_{d+2}=x_1+x_2\mbox{ and }x_d+x_{d+3}=0.$

\bigskip

(4) $\alpha=1$. The primitive relations are $x_1+\cdots+x_{d-1}=x_d,\ x_2+\cdots+x_{d-1}+x_{d+3}=0,\ x_d+x_{d+3}=x_1$ and
\begin{center}
\begin{tabular}{|c||c|c|c|} \hline
Case & 1 & 2 & 3 \\ \hline\hline
$x_1+x_{d+1}+x_{d+2}=$ & $x_{d+3}$ & $2x_{d+3}$ & $x_2+x_{d+3}$ \\ \hline
$x_d+x_{d+1}+x_{d+2}=$ & $0$ & $x_{d+3}$ & $x_2$ \\ \hline
\end{tabular}
\end{center}

\bigskip

(IIIa) \underline{$\rho(X)=4$ and $X$ is an $S_7$-bundle over ${\bf P}^{d-2}$.} Let $\G(\Sigma)=\{x_1,\ldots,x_{d+4}\}$. The primitive relations of $\Sigma$ are $x_d+x_{d+2}=x_{d+1},\ x_{d}+x_{d+3}=x_{d+4},\ x_{d+1}+x_{d+3}=0,\ x_{d+1}+x_{d+4}=x_{d},\ x_{d+2}+x_{d+4}=0$ and
\begin{center}
\begin{tabular}{|c||c|c|c|} \hline
Case & 1 & 2 & 3 \\ \hline\hline
$x_1+\cdots+x_{d-1}=$ & $\alpha x_d$ & $\alpha x_{d+1}$ & $\alpha x_{d+2}$ \\ \hline
\end{tabular}
\end{center}

\bigskip

(IIIb) \underline{$\rho(X)=4$ and $X$ is a toric bundle over ${\bf P}^1$.} $X$ is an $X'$-bundle over ${\bf P}^1$, where $X'$ is a toric Fano $(d-1)$-fold equipped with an extremal contraction which contracts a divisor to a point. Let $\G(\Sigma)=\{x_1,\ldots,x_{d+4}\}$. The primitive relations of $\Sigma$ are $x_1+\cdots+x_{d-1}=\alpha x_{d},\ x_2+\cdots+x_{d-1}+x_{d+4}=(\alpha-1)x_{d},\ x_{d}+x_{d+3}=0,\ x_d+x_{d+4}=x_1,\ x_1+x_{d+3}=x_{d+4}$ and
\begin{center}
\begin{tabular}{|c||c|c|c|c|c|c|} \hline
Case & 1 & 2 & 3 & 4 & 5 & 6 \\ \hline\hline
$x_{d+1}+x_{d+2}=$ & $0$ & $x_1$ & $x_2$ & $x_d$ & $x_{d+3}$ & $x_{d+4}$ \\ \hline
\end{tabular}
\end{center}

\bigskip

(IV) \underline{$\rho(X)=5$.} $X$ is an ${S_6}$-bundle over ${\bf P}^{d-2}$. The primitive relations of $\Sigma$ are $x_1+\cdots+x_{d-1}={\alpha}x_d,\ x_d+x_{d+2}=x_{d+1},\ x_d+x_{d+3}=0,\ x_d+x_{d+4}=x_{d+5},\ x_{d+1}+x_{d+3}=x_{d+2},\ x_{d+1}+x_{d+4}=0,\ x_{d+1}+x_{d+5}=x_d,\ x_{d+2}+x_{d+4}=x_{d+3},\ x_{d+2}+x_{d+5}=0$ and $x_{d+3}+x_{d+5}=x_{d+4}$, where $\G(\Sigma)=\{x_1,\ldots,x_{d+5}\}$.

\begin{Rem}

{\rm
If $d\leq 4$, there exist toric Fano $d$-folds equipped with extremal contractions which contract divisors to curves which are not contained in this list $($see Batyrev \cite{batyrev4}, Sato \cite{sato1} and Watanabe-Watanabe \cite{watanabe1}$)$.
}

\end{Rem}

%%%%%%%%%%%%%%%%%%%%%%%%%%%%%%%%%%%%%%%%%%%%%
\section{Projective toric varieties whose toric blow-up along a curve is Fano}\label{table2}
%%%%%%%%%%%%%%%%%%%%%%%%%%%%%%%%%%%%%%%%%%%%%

\hspace{.5cm} By blowing-down the toric Fano $d$-folds in Section \ref{table1} such that $\alpha=1$, we can obtain the classification of smooth projective toric $d$-folds which can be equivariantly blown-up to Fano along curves. We assume $d\geq 5$.

\bigskip

(I) \underline{$\rho(X)=1$.} $X$ is ${\bf P}^d$. The primitive relation is
$$x_1+\cdots+x_{d+1}=0,$$
where $\G(\Sigma)=\{x_1,\ldots,x_{d+1}\}$.

\bigskip

(IIa) \underline{$X$ is either a ${\bf P}^1$-bundle over ${\bf P}^{d-1}$ or a ${\bf P}^{d-1}$-bundle over ${\bf P}^{1}$.} Let $\G(\Sigma)=\{x_1,$ $\ldots,x_{d+2}\}$. The primitive relations are as follows:
\begin{center}
\begin{tabular}{|c||c|c|c|c|} \hline
Case & 1 & 2 & 3 & 4 \\ \hline\hline
$x_1+\cdots+x_{d}=$ & $0$ & $x_{d+1}$ & $0$ & $0$ \\ \hline
$x_{d+1}+x_{d+2}=$ & $0$ & $0$ & $x_1+\cdots+x_{d-1}$ & $x_1$ \\ \hline
\end{tabular}
\end{center}

\bigskip

(IIb) \underline{$\rho(X)=2$ and $X$ is a ${\bf P}^{d-2}$-bundle over ${\bf P}^2$.} Let $\G(\Sigma)=\{x_1,\ldots,x_{d+2}\}$. The primitive relations are $x_1+\cdots+x_{d-1}=0$ and
\begin{center}
\begin{tabular}{|c||c|c|c|} \hline
Case & 1 & 2 & 3 \\ \hline\hline
$x_d+x_{d+1}+x_{d+2}=$ & $x_1$ & $2x_1$ & $x_1+x_2$ \\ \hline
\end{tabular}
\end{center}

\bigskip

(IIIa) \underline{$\rho(X)=3$ and $\Sigma$ is not a splitting fan.} Let $\G(\Sigma)=\{x_1,\ldots,x_{d+3}\}$. The primitive relations of $\Sigma$ are as follows:
\begin{center}
\begin{tabular}{|c||c|c|c|} \hline
Case & 1 & 2 & 3 \\ \hline\hline
$x_{d+2}+x_{d+3}=$ & $x_1+\cdots+x_{d-1}$ & $x_1+\cdots+x_{d-1}$ & $0$ \\ \hline
$x_1+\cdots+x_{d}=$ & $x_{d+2}$ & $0$ & $x_{d+2}$ \\ \hline
$x_1+\cdots+x_{d-1}+x_{d+1}=$ & $x_{d+3}$ & $x_{d+3}$ & $0$ \\ \hline
$x_{d}+x_{d+3}=$ & $0$ & $x_{d+1}$ & $x_{d+1}$ \\ \hline
$x_{d+1}+x_{d+2}=$ & $0$ & $0$ & $x_{d}$ \\ \hline
\end{tabular}
\end{center}

\bigskip

(IIIb) \underline{$\rho(X)=3$ and $\Sigma$ is a splitting fan.} Let $\G(\Sigma)=\{x_1,\ldots,x_{d+3}\}$. The primitive relations of $\Sigma$ are $x_2+\cdots+x_{d}=0,\ x_1+x_{d+3}=x_d$ and
\begin{center}
\begin{tabular}{|c||c|c|c|c|c|} \hline
Case & 1 & 2 & 3 & 4 & 5 \\ \hline\hline
$x_{d+1}+x_{d+2}=$ & $0$ & $x_1$ & $x_2$ & $x_1+\cdots+x_{d-1}$ & $x_d$ \\ \hline
\end{tabular}
\end{center}

\bigskip

(IV) \underline{$\rho(X)=4$.} The primitive relations of $\Sigma$ are $x_1+\cdots+x_{d-1}+x_{d+2}=x_{d+1},\ x_1+\cdots+x_{d-1}+x_{d+3}=0,\ x_1+\cdots+x_{d-1}+x_{d+4}=x_d,\ x_{d+1}+x_{d+3}=x_{d+2},\ x_{d+1}+x_{d+4}=0,\ x_d+x_{d+1}=x_1+\cdots+x_{d-1},\ x_{d+2}+x_{d+4}=x_{d+3},\ x_d+x_{d+2}=0$ and $x_d+x_{d+3}=x_{d+4}$, where $\G(\Sigma)=\{x_1,\ldots,x_{d+4}\}$.

\begin{Rem}

{\rm
If $d\leq 4$, there exist smooth projective toric $d$-folds which can be equivariantly blown-up to Fano along curves which are not contained in this list $($see Batyrev \cite{batyrev4}, Sato \cite{sato1} and Watanabe-Watanabe \cite{watanabe1}$)$.
}

\end{Rem}

\bigskip

\begin{flushleft}
\begin{sc}
Department of Mathematics \\
Tokyo Institute of Technology \\
Oh-Okayama, Meguro, Tokyo \\
Japan
\end{sc}

\medskip
{\it E-mail address}: $\mathtt{hirosato@math.titech.ac.jp}$
\end{flushleft}


\begin{thebibliography}{99}



\bibitem{batyrev3} V. V. Batyrev, On the classification of smooth projective toric varieties, Tohoku Math. J. 43 (1991), 569--585.


\bibitem{batyrev4} V. V. Batyrev, On the classification of toric Fano 4-folds, Algebraic geometry, 9, J. Math. Sci. (New York) 94 (1999), 1021--1050.

\bibitem{bonavero1} L. Bonavero, Toric varieties whose blow-up at a point is Fano, math.AG/0012229.

\bibitem{bonavero2} L. Bonavero, F. Campana, J. A. Wi\'sniewski, Vari\'et\'es projectives complexes dont l'\'eclat\'ee en un point est de Fano, math.AG/0106047.



\bibitem{casagrande1} C. Casagrande, Contractible classes in toric varieties, math.AG/0111332.

\bibitem{casagrande2} C. Casagrande, Toric Fano varieties and birational morphisms, math.AG/0112007.




\bibitem{fulton1} W. Fulton, Introduction to Toric Varieties, Ann. of Math. Studies 131, Princeton Univ. Press, Princeton, NJ, 1993.

\bibitem{klein1} P. Kleinschmidt, A classification of toric varieties with few generators, Aequationes Math. 35 (1988), 254--266.






\bibitem{oda2} T. Oda, Convex Bodies and Algebraic Geometry---An introduction to the theory of toric varieties, Ergeb. Math. Grenzgeb. (3), Vol. 15, Springer-Verlag, Berlin, Heidelberg, New York, London, Paris, Tokyo, 1988.


\bibitem{reid1} M. Reid, Decomposition of toric morphisms, in Arithmetic and Geometry, papers dedicated to I. R. Shafarevich on the occasion of his 60th birthday (M. Artin and J. Tate, eds.), vol. II, Geometry, Progress in Math. 36, Birkh\"{a}user, Boston, Basel, Stuttgart, 1983, 395--418.

\bibitem{sato1} H. Sato, Toward the classification of higher-dimensional toric Fano varieties, Tohoku Math. J. 52 (2000), 383--413.


\bibitem{watanabe1} K. Watanabe and M. Watanabe, The classification of Fano 3-folds with torus embeddings, Tokyo J. Math. 5 (1982), 37--48.



\end{thebibliography}
\end{document}